\documentclass{eajam}
\renewcommand\thisnumber{x}
\renewcommand\thisyear {200x}
\renewcommand\thismonth{xxx}
\renewcommand\thisvolume{xx}
\usepackage{charter}
\usepackage[charter]{mathdesign}

\usepackage{cleveref}
\usepackage{amsmath,amssymb,amsthm,mathrsfs}
\usepackage[caption=false,font=normalsize,labelfont=sf,textfont=sf]{subfig}
\usepackage{booktabs,makecell,multirow}
\usepackage{dsfont}

\crefname{equation}{Eq.}{}
\crefname{section}{Sect.}{}
\crefname{figure}{Fig.}{}
\crefname{table}{Table.}{}

\newcommand{\Nm}[1]{\left\lVert {#1} \right\rVert}

\newtheorem{discretization}{Discretization}[section]
\begin{document}

\markboth{Maoyuan Xu and Xiaoping Xie}{An efficient feature-preserving PDE for image denoising}
\title{An efficient feature-preserving PDE algorithm for image denoising based on a spatial-fractional anisotropic diffusion equation}

\author[Maoyuan Xu and Xiaoping Xie]{Maoyuan Xu\affil{1} and Xiaoping Xie\affil{1}\comma\corrauth}
\address{\affilnum{1}\ School of Mathematics, Sichuan University, Chengdu 610065, China
}
\emails{{\tt 663521583@qq.com} (M. Xu), {\tt xpxie@scu.edu.cn} (X. Xie)}

\begin{abstract}
How to effectively remove the noise while preserving the image structure features is a challenging issue in the field of image denoising. In recent years,
fractional PDE based methods have attracted more and more research efforts due to the ability to balance the noise removal and the preservation of image edges and textures.
Among the existing fractional PDE algorithms, there are only a few using spatial fractional order derivatives, and  all the fractional derivatives involved are one-sided derivatives. In this paper, an efficient  feature-preserving fractional PDE algorithm is proposed for image denoising based on a  nonlinear spatial-fractional anisotropic diffusion equation.  Two-sided Gr\"umwald-Letnikov fractional derivatives were used in the PDE model which are suitable to depict the local self-similarity of images. The  Short Memory Principle is employed to simplify the approximation scheme. Experimental results show that the proposed method is of a satisfactory performance, i.e. it keeps a remarkable balance between noise removal and feature preserving, and has an extremely high structural retention property.
\end{abstract}



\maketitle

\section{Introduction}
With the development of modern society, the form and quantity of information are growing rapidly, and most of the information are images.
In the process of acquisition, transmission and storage, the images inevitably  suffer from noise interference. Such interference will blur the images and cause the difficulty to identify interested targets, thus seriously affecting the subsequent processing and analysis of the images.
Nowadays image denoising has become a particularly important field in  image processing, and its  challenging aspect is how to effectively remove the noise while preserving the image structure features.

There are many types of algorithms for image denoising, such as convolutional filter based methods (\cite{pu2009fractional,eng2001noise}), wavelet analysis based methods (\cite{quiroga2003single,chen2005image}), machine learning based methods (\cite{xie2012image,zhang2017beyond}), and variation or PDE based methods. Generally speaking, a conventional convolutional filter based method is easy to conduct, but often has ordinary performance; a wavelet analysis based method is of better performance when selecting an appropriate threshold, but the Gibbs phenomenon will appear at discontinuous points of the signal; when using machine learning based method, there is a great demand of high quality training samples, and designing an appropriate network architecture is also a skilful work.

Compared with other methods, the variation or PDE (partial differential equation) based methods are based on rigorous mathematics and are easier to explain, though it is technical to design effective formulas for the methods. In a variation based method, according to some mathematics analysis of the images, an energy functional is given firstly to obtain an optimization problem which restores the observation images to the original clean ones, then a numerical algorithm is adopted to solve the optimization problem so as to acquire the processed images; see, e.g. \cite{sapiro2006geometric,chan2000high,han2019tensor,vogel1996iterative,chen2010adaptive,getreuer2012rudin,chan2001total,scherzer2009variational,rudin1992nonlinear,pu2014fractional1,pu2014fractional2}.

In a PDE based method, the PDE employed for image denoising is usually constructed directly on the basis of the mathematical properties of the diffusion and the structural characteristics of the image. For instance, Perona and Malik \cite{perona1990scale} introduced an algorithm based on a diffusion process (later called the PM model), where the diffusion coefficient is chosen to vary spatially so as to encourage intra-region smoothing in preference to inter-region smoothing.
Experimental results showed that their method has remarkable performance on edge preserving when denoising; nevertheless, the PDE of their method is a second-order partial differential equation, hence the processed images suffer from the staircase effect (the piecewise constant effect)  which makes the image artificial. We note that  there also have developed many other methods based on second-order PDEs, e.g. \cite{barbu2015nonlinear,rafsanjani2017adaptive,tang2001color,gilboa2004image,sapiro2006geometric}, but the staircase effect still happens in these methods. In order to eliminate the staircase effect, some fourth-order PDEs based methods were proposed   for image denoising (\cite{you2000fourth,hajiaboli2011anisotropic, barbu2016nonlinear,deng2019hessian}); yet these methods smooth the image too much to keep discontinuity, and  the speckle effect appears.

In recent years, fractional order PDE based methods have attracted more and more research interest due to their ability to balance the noise removal and the preservation of image edges and textures; see, e.g. \cite{cuesta2003image, yu2018image, abirami2018fractional, janev2011fully, liao2013low, cuesta2012image, zhang2015spatial}. Numerical experiments reveal that these methods are able to alleviate both the staircase and speckle effects and generate processed images of higher quality than the integer order PDE based methods. It should be mentioned  that  fractional order derivatives were also used in some variation based methods (cf. \cite{bai2007fractional,bai2018image,pu2014fractional1,pu2014fractional2,dong2013variational, dong2016fractional, zhang2015total, jia2017new, litvinov2011modified}), where the corresponding Euler equations are fractional order PDEs.

So far, among the existing fractional PDE based methods, there are only a few using spatial fractional order derivatives. In \cite{liao2013low}, a spatial fractional order anisotropic diffusion equation like
\begin{equation*}
  \frac{\partial u}{\partial t} = \text{div}^{\alpha}\Big[c\Big(\Nm{G_{\sigma_1}\cdot\nabla^{\alpha}u}\Big)\nabla^{\alpha}u\Big]
\end{equation*}
was applied for low-dose computed tomography imaging, where  the fractional divergence and gradient operators, $\text{div}^{\alpha}$ and $\nabla^{\alpha}$, are defined with  Gr\"umwald-Letnikov fractional derivatives of order $\alpha\in [0.5,1.5]$.
In \cite{zhang2015spatial} the following spatial fractional telegraph equation was used for image structure preserving denoising:
\begin{equation*}
   \frac{\partial ^2u}{\partial t^2} + \lambda  \frac{\partial u}{\partial t} - \text{div}^{\alpha}\Big(g\Big(|\nabla^{\alpha}G_{\sigma}\ast u|\Big)\nabla^{\alpha}u\Big) = 0,
\end{equation*}
where $\text{div}^{\alpha}$ and $\nabla^{\alpha}$ are given in terms of Riemann-Liouville fractional derivatives of order $\alpha\in [1,2]$.
On the basis of the Rudin-Osher-Fatemi model (ROF or TV model) \cite{rudin1992nonlinear}, a spatial fractional isotropic diffusion equation of the form
\begin{equation*}
  \frac{\partial u}{\partial t} = \text{div}^{\alpha} u
\end{equation*}
was employed in \cite{abirami2018fractional} for image denoising, with $\text{div}^{\alpha}$ being defined by Gr\"umwald-Letnikov fractional derivatives of order $\alpha\in (1,2)$.
It should be pointed out that all the spatial fractional order derivatives used in \cite{abirami2018fractional, liao2013low, zhang2015spatial} are one-sided derivatives which, as shown in numerical experiments, usually result in underutilizing of image pixels.

In this paper, we shall apply the following spatial fractional order anisotropic diffusion equation for image denoising:
\begin{equation*}
 \frac{\partial u}{\partial t} = -\text{div}^{\alpha}(g(|\nabla^{\beta}u|)\nabla^{\alpha}u),
\end{equation*}
where $1.25 < \alpha < 1.75, \ 1 < \beta < 2,$ and  $\text{div}^{\alpha}$, $\nabla^{\alpha}$  and $\nabla^{\beta}$ are defined by   two-sided Gr\"umwald-Letnikov fractional derivatives of orders $\alpha\in (1.25, 1.75)$ and  $\beta\in (1,2)$; see \cref{section:preliminaries} and \cref{section:Proposed diffusion equation} for details. We note that, from the mathematical point of view, the use of a two-sided or symmetric-form fractional derivative involves more neighbourhood information and is more suitable to depict the local self-similarity of images than that of a one-sided derivative, and thus is expectable to avoid underutilizing of image pixels in image denoising.

The rest of the paper is organized as follows. \cref{section:preliminaries} introduces the definitions of Gr\"umwald-Letnikov fractional derivative and fractional divergence and gradient operators, and discusses the properties of fractional order derivative. \cref{section:Proposed diffusion equation} proposes the spatial-fractional anisotropic diffusion model and gives the numerical discretization scheme of the fractional PDE model. In \cref{section:Experimental results}, experimental results are depicted to test the performance of the proposed method. Finally, \cref{section:Conclusion} presents some concluding remarks.

\section{Preliminaries}\label{section:preliminaries}
\subsection{Gr\"umwald-Letnikov fractional derivatives}
There are many definitions of fractional order derivatives, among which Gr\"umwald-Letnikov, Riemann-Liouville and Caputo definitions are the most common ones. In our model we will use  the Gr\"umwald-Letnikov definition.

Let $n>0$ and $ L>0$ be a integer and a constant, respectively. For a function  $f(x)\in C^n[0,L]$, it is well known that its derivative of order $n$ can be expressed as the limit of backward differences difference quotient of order $n$  with step size $h$, i.e.
\begin{align}\label{derivative-n}
  \frac{d^nf(x)}{dx^n} =& \lim_{h\rightarrow 0}\frac{1}{h^n}\sum_{k=0}^{n}(-1)^k\binom{n}{k}f(x-kh) \nonumber\\
                       =& \lim_{h\rightarrow 0}\frac{1}{h^n}\sum_{k=0}^{\infty}(-1)^k\binom{n}{k}f(x-kh) , \ x\in (0,L].
\end{align}
where $\binom{r}{k}$, for any real number $r$ and positive integer $k\geq 0$, denote the binomial coefficients defined by
\begin{eqnarray*}
 \binom{r}{k}:=\left\{\begin{array}{ll} \frac{r(r-1)(r-2)\cdots (r-k+1)} {k!} =\frac{\Gamma(r+1)}{\Gamma(k+1)\Gamma(r-k+1)}&\text{if } k\leq r,\\
 0&\text {if } k>r,
 \end{array}
 \right.
\end{eqnarray*}
and $\Gamma(\cdot)$ is the Gamma function with $\Gamma(z)=\int_0^\infty t^{z-1}e^{-t}dt$ if $z>0$, and $\Gamma(z)=z^{-1}\Gamma(z+1)$ if $z<0$.

By extending the formulation \eqref{derivative-n} from  the integer order $n$ to a real number $\alpha\in (n-1,n], $ we can obtain the definition of Gr\"umwald-Letnikov fractional  derivative.

\begin{definition}\label{eq:definition2} The left-sided Gr\"{u}mwald-Letnikov fractional derivative of order $\alpha$ of the function $f$ is defined by
\begin{equation*}
  {}_0D^{\alpha}_xf(x): = \lim_{h\rightarrow 0}\frac{1}{h^\alpha}\sum_{k=0}^{[x/h]}\omega^{(\alpha)}_kf(x-kh), \ x\in (0,L],
\end{equation*}
where
\begin{equation*}
\omega_k^{(\alpha)}=(-1)^k\binom{\alpha}{k}=\frac{(-1)^k\Gamma(\alpha+1)}{\Gamma(k+1)\Gamma(\alpha-k+1)}=\frac{\Gamma(k-\alpha)}{\Gamma(k+1)\Gamma(-\alpha)}
\end{equation*}
are called the Gr\"{u}mwald-Letnikov coefficients. The right-sided Gr\"{u}mwald-Letnikov fractional derivative of order $\alpha$ of the function $f$ is defined by
\begin{equation*}
  {}_xD^{\alpha}_L f(x) := \lim_{h\rightarrow 0}\frac{1}{h^\alpha}\sum_{k=0}^{[(L-x)/h]}\omega^{(\alpha)}_kf(x+kh) , \ x\in [0,L).
\end{equation*}
\end{definition}

Taking the local self-similarity of images into account, we also need to introduce the following symmetric two-sided Gr\"{u}mwald-Letnikov fractional derivative.

\begin{definition}\label{eq:definition-two-sided} The two-sided Gr\"{u}mwald-Letnikov fractional derivative of order $\alpha$ of the function $f$ is defined by
\begin{align*}
  D^{\alpha}_xf: =& \frac{1}{2}({}_0D^{\alpha}_x+{}_xD^{\alpha}_L)f, \  x\in (0,L).
\end{align*}
\end{definition}

Based on \cref{eq:definition-two-sided}, we introduce the following definitions of the two-sided fractional gradient operator $\nabla^{\alpha}$ and divergence operator $\text{div}^{\alpha}$.
\begin{definition}\label{eq:definition-div-grad}
For a real scalar function $u(x,y)$ and a real vector function $\vec{f} = \left(f_1(x,y),f_2(x,y)\right)$ with $x,y\in (0,L)$, define
\begin{align*}
   \nabla^{\alpha}u: =& (D^{\alpha}_xu, D^{\alpha}_yu), \\
      \text{div}^{\alpha}\vec{f} :=&\nabla^{\alpha}\cdot \vec{f} = D^{\alpha}_x f_1+ D^{\alpha}_y f_2.
\end{align*}
\end{definition}

In addition, we introduce several norms as follows.
\begin{definition}\label{def:norm1}
For function $u(x,y)$, matrix $A\in \mathds{R}^{M\times N}$ and vector $\upsilon\in\mathds{R}^{M\times1}$, define
\begin{align*}
  |\nabla^{\alpha}u| :=& \sqrt{{(D^{\alpha}_x u)}^2 + {(D^{\alpha}_y u)}^2}, \\
  \Nm{A}_{\infty} :=& \max_{1\leq i\leq M}\sum^{N}_{j=1}|a_{i,j}|, \\
  \Nm{\upsilon}_{\infty}: =& \max_{1\leq i\leq M}|\upsilon_i|.
\end{align*}
\end{definition}

\subsection{Performance of fractional derivatives}\label{subsection:Performanceoffractionalorderderivatives}
As mentioned before, the fractional PDE based methods are more capable of  balancing the noise removal and the preservation of image edges and textures in image denoising than the integer order PDE based methods. The main reason for this lies in that the fractional differential is able to nonlinearly preserve the low frequency information, which contains abundant surface texture details, and to enhance the high frequency information of images (cf. \cite{pu2009fractional}), while the integer order derivatives behave not well in dealing with the low frequency information. In fact, take the left-sided Gr\"{u}mwald-Letnikov fractional derivative $ {}_0D^{\alpha}_xf(x)$, $n-1<\alpha\leq n$, as an example and let $\mathcal{F}$ be the Fourier transform, then
\begin{equation}\label{eq:Fractionalproperty1}
   \mathcal{F}[{}_0D^{\alpha}_xf(x)] (\omega) = (i\omega)^{\alpha} \mathcal{F}[f(x)](\omega) - \sum^{n-1}_{k=0}(i\omega)^k \frac{d^{\alpha-k-1}}{dx^{\alpha-k-1}}f(0),
\end{equation}
where for $\gamma=\alpha-k-1$,
\begin{equation*}
 \frac{d^{\gamma}}{dx^{\gamma}}f(x) := \frac{1}{\Gamma(-\gamma)}\int_0^x (x-\tau)^{-\gamma-1}f(\tau)d\tau.
\end{equation*}
Note that the term
\begin{equation*}
\sum^{n-1}_{k=0}(i\omega)^k\frac{d^{\alpha-k-1}}{dx^{\alpha-k-1}}f(0)
\end{equation*}
in \cref{eq:Fractionalproperty1} does not alter the waveform of frequency spectrum of images, but the factor $(i\omega)^{\alpha} $ does.

\begin{figure}[htbp!]
  \centering
  \includegraphics[width = 3in]{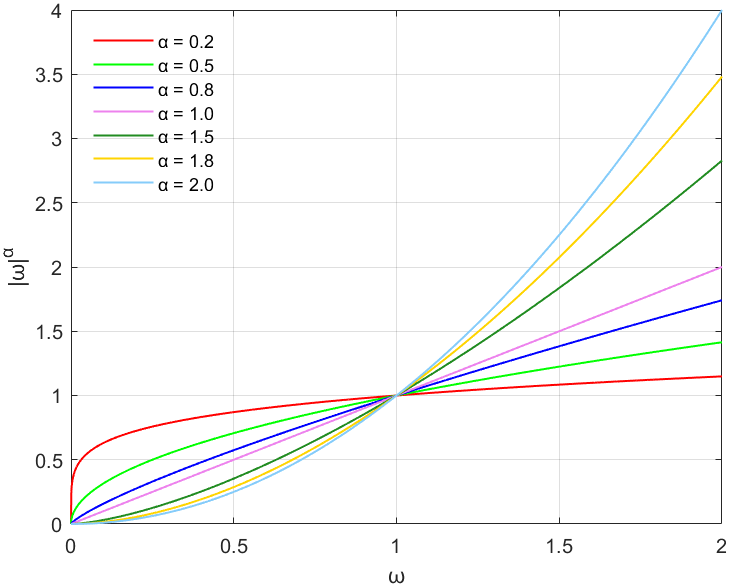}\\
  \caption{Frequency response of fractional order derivative.}\label{fig:FrequencyResponse1}
\end{figure}

\cref{fig:FrequencyResponse1} shows the amplitude frequency characteristic curves of fractional order derivatives with different fractional orders. We can see that in the low frequency region of $0<\omega<1$, the smaller the the order $\alpha$ of fractional derivative is, the stronger the capacity to preserve the magnitude will be; and in the high frequency region of $\omega>1$, the fractional  derivative of order $\alpha>1$ is superior to the  first order derivative in enhancing  the high frequency components.

\begin{figure}[htbp!]
  \centering
  \subfloat[Original image $u$ ]{\includegraphics[width= 1.3in]{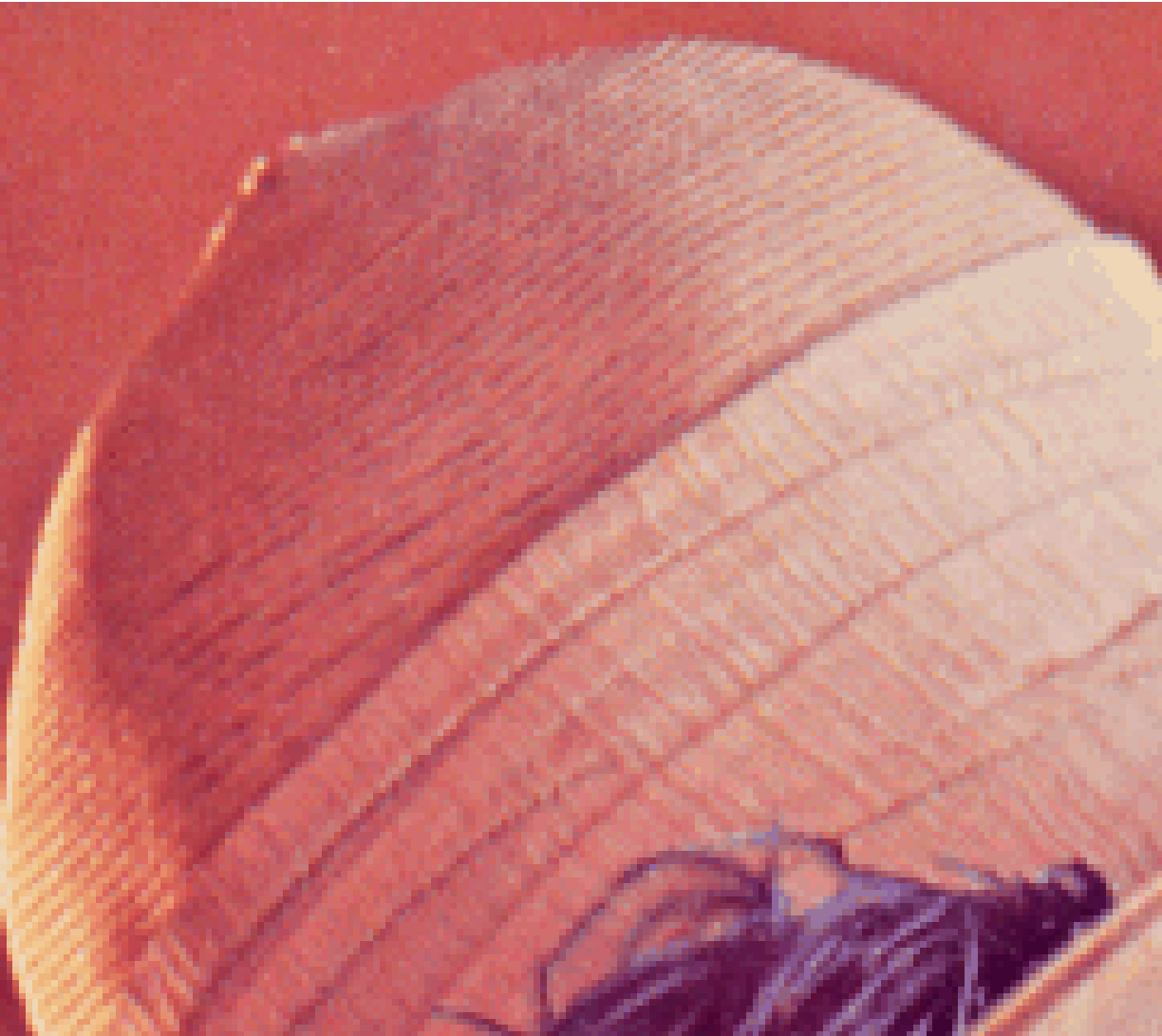}}\quad\quad
  \subfloat[$|\nabla u|$ image]{\includegraphics[width= 1.3in]{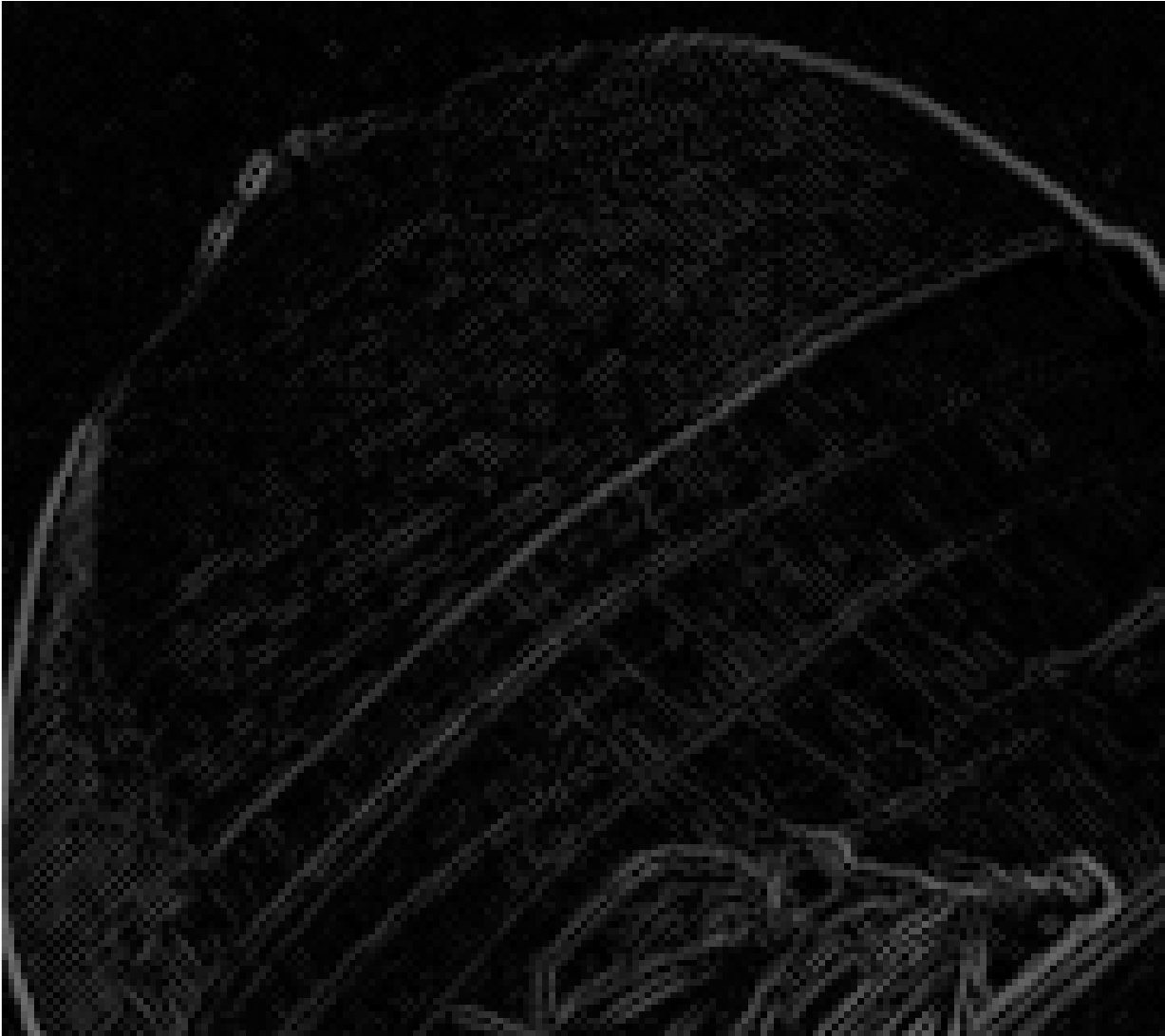}}\quad\quad
  \subfloat[$|\nabla^{1.5} u|$ image]{\includegraphics[width= 1.3in]{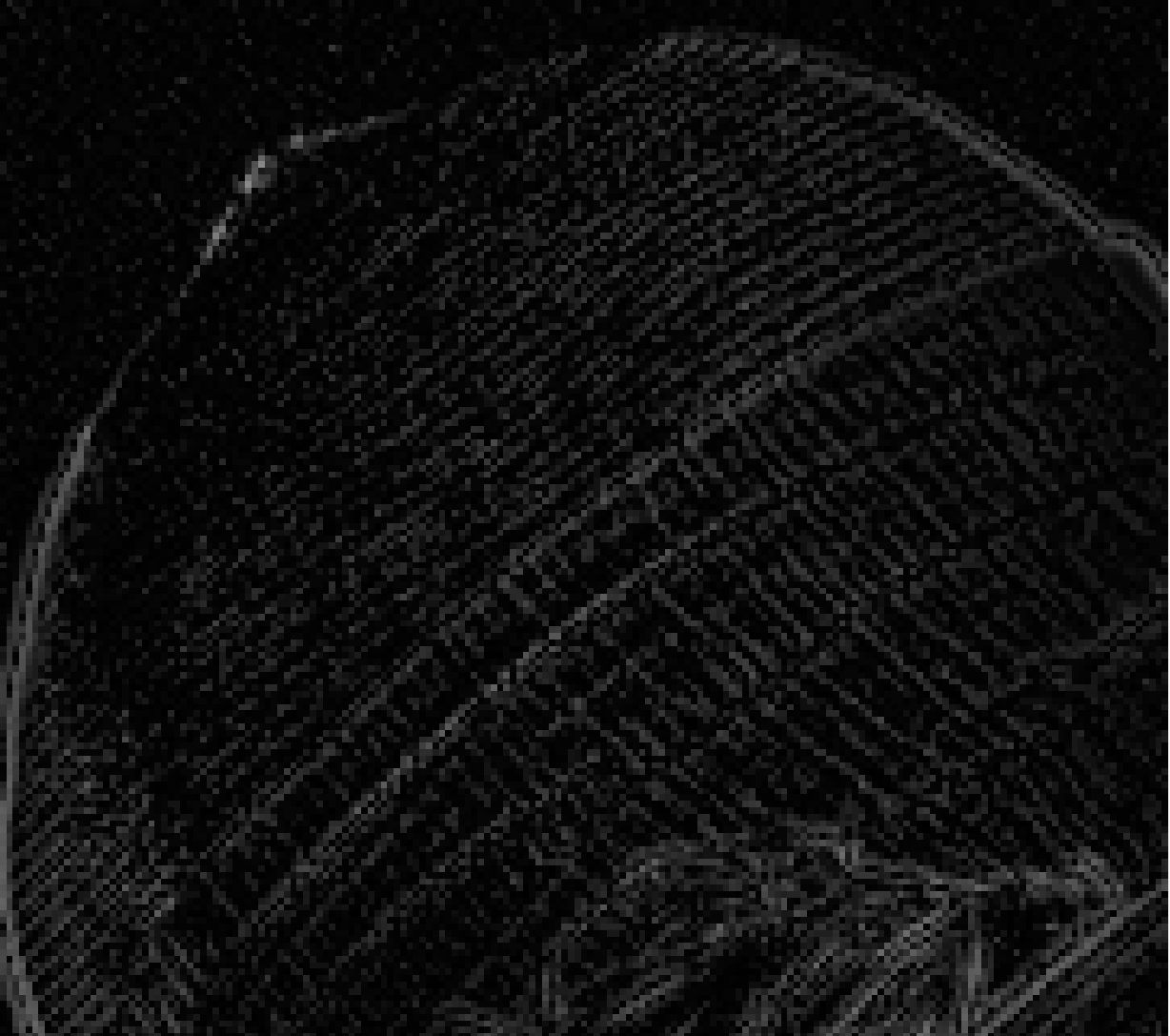}}\\
  \caption{Comparison between $|\nabla u|$ and $|\nabla^{1.5} u|$.}\label{fig:comparison1}
\end{figure}

On the other hand, the edges and textures of image are both deemed as the most important features of the image, since they contain a large quantity of components with rich high frequency. Therefore, it is ideal to remove the noise without losing any features; Nevertheless, traditional methods such as the Gaussian filter and the Median filter have poor performance at feature preserving. In the classic PM model, $|\nabla u|$ is employed to detect the edges and textures of the image, but the first order based differential operator is not able to gather ample features. As shown in \cref{fig:comparison1}, the feature detecting capability of the two-sided fractional order derivative based gradient term $|\nabla^{\alpha} u|$ with $\alpha=1.5$ is more superior than that of $|\nabla u|$.

\section{Spatial-fractional diffusion model and discretization scheme}\label{section:Proposed diffusion equation}
\subsection{Spatial-fractional anisotropic diffusion model}
\begin{figure}[htbp!]
  \centering
  \includegraphics[width= 3in]{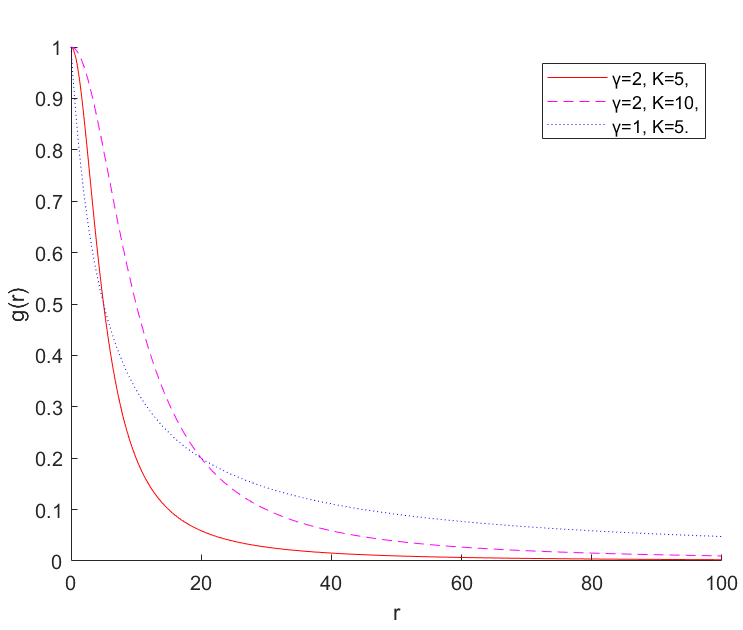}\\
  \caption{Graph of $g(r)$ in \cref{eq:g1} and \cref{eq:g2}.}\label{fig:g1}
\end{figure}

Let $\Omega=[0,L]\times [0,L]$ be the compact support of the image domain. Inspired by the analysis in \cref{subsection:Performanceoffractionalorderderivatives}, we propose the following spatial-fractional order two-sided derivative based diffusion model:
\begin{equation}\label{eq:origin}
\left\{
\begin{aligned}
  & \frac{\partial u}{\partial t} = -\text{div}^{\alpha}(g(|\nabla^{\beta}u|)\nabla^{\alpha}u), &in ~\Omega\times[0,T], \\
  & \frac{\partial u}{\partial \vec{n}} = 0, &on ~\partial\Omega\times[0,T], \\
  & u(x,y,0) = u_0(x,y), &\forall  (x,y)\in\Omega.
\end{aligned}
\right.
\end{equation}
Here $ 1.25 < \alpha < 1.75, ~1 < \beta < 2,$  $\vec{n}$ is the unit outer normal vector along $\partial\Omega$, and $u_0$ is the initial image to be processed. $g(r)$ is an edge-stopping function used to control the rate of diffusion to retain more image features while diffusing, which is assumed to satisfy the following conditions (cf. \cite{perona1990scale}):
\begin{equation*}
\left\{
\begin{aligned}
  (a)&  ~g(r)~\text{is a non-increasing function}; \\
  (b)&  ~\lim_{r\rightarrow 0} ~g(r) = 1; \\
  (c)&  ~\lim_{r\rightarrow\infty} g(r) = 0.
\end{aligned}
\right.
\end{equation*}
We note that in \cite{perona1990scale} two examples of $g(r)$ are provided:
\begin{align}
     g(r) & = \frac{1}{1+(\frac{r}{K})^{\gamma}},\label{eq:g1} \\
     g(r) & = \exp\Big(-\Big(\frac{r}{K}\Big)^{\gamma}\Big) \label{eq:g2},
\end{align}
where $K>0$ is a threshold to control the attenuation rate of $g(r)$, and $\gamma = 1,2$. \cref{fig:g1} shows the graph of $g(r)$ in \cref{eq:g1} and \cref{eq:g2} with different $\gamma$ and $K$.

\subsection{Numerical discretization scheme}\label{section:Numerical discretization scheme}
We first introduce approximation formulas for the left-sided and right-sided Gr\"{u}mwald-Letnikov fractional derivatives. In \cite{oldham;1974}, based on the relation
\begin{equation}\label{eq:approximation1}
  {}_0D^{\alpha}_x f(x) = \lim_{h\rightarrow 0}\frac{h^{-\alpha}}{\Gamma(-\alpha)}\sum_{k=0}^{[x/h+\alpha/2]}\frac{\Gamma(k-\alpha)}{\Gamma(k+1)}f(x-(k-\frac{1}{2}\alpha)h)
\end{equation}
and the three-points Lagrange interpolation formula
\begin{equation*}
\begin{split}
  f(x-(k-\frac{1}{2}\alpha)h)\approx & (\frac{\alpha}{4}+\frac{\alpha^2}{8})f(x-(k-1)h) \\
                                     & +(1-\frac{\alpha^2}{4})f(x-kh) \\
                                     & +(\frac{\alpha^2}{8}-\frac{\alpha}{4})f(x-(k+1)h)
\end{split}
\end{equation*}
for any $0<x \leq L$ and sufficiently small $h$, the following approximation formula of second order accuracy, called G2 scheme, was proposed to compute the left-sided Gr\"{u}mwald-Letnikov fractional derivative:
\begin{equation}\label{eq:fracderivative1}
\begin{split}
    {}_0D^{\alpha}_x f(x) \approx & h^{-\alpha}\sum^{n_1-1}_{k=0}\omega_{k}^{(\alpha)}\Big\{f(x-kh) \\
                              & +\frac{1}{4}\alpha\Big(f(x-(k-1)h)-f(x-(k+1)h)\Big) \\
                              & +\frac{1}{8}\alpha^2\Big(f(x-(k-1)h)-2f(x-kh) \\
                              &\quad\quad\quad +f(x-(k+1)h)\Big)\Big\},
\end{split}
\end{equation}
where $n_1=[x/h+\alpha/2]+1$, and $\omega_k^{(\alpha)}$ is given in \cref{eq:definition2}.

Similarly, from the relation
\begin{equation}\label{eq:approximation1}
 {}_0D^{\alpha}_x f(x) = \lim_{h\rightarrow 0}\frac{h^{-\alpha}}{\Gamma(-\alpha)}\sum_{k=0}^{[(L-x)/h+\alpha/2]}\frac{\Gamma(k-\alpha)}{\Gamma(k+1)}f(x+(k-\frac{1}{2}\alpha)h),
\end{equation}
the G2 scheme for the right-sided Gr\"{u}mwald-Letnikov fractional derivative is given by
\begin{equation}\label{eq:fracderivative1-R}
\begin{split}
 {}_xD^{\alpha}_L f(x) \approx & h^{-\alpha}\sum^{n_2-1}_{k=0}\omega_{k}^{(\alpha)}\Big\{f(x+kh) \\
                              & +\frac{1}{4}\alpha\Big(f(x+(k-1)h)-f(x+(k+1)h)\Big) \\
                              & +\frac{1}{8}\alpha^2\Big(f(x+(k-1)h)-2f(x+kh) \\
                              &\quad\quad\quad +f(x+(k+1)h)\Big)\Big\}
 \end{split}
\end{equation}
for any $0\leq x < L$ and sufficiently small $h$, where $n_2=[(L-x)/h+\alpha/2]+1$.

The Short Memory Principle \cite{podlubny1997numerical} is to be employed to simplify the approximation formulas \cref{eq:fracderivative1} and \cref{eq:fracderivative1-R}.
In fact, when $x\gg0$, there are lots of summation terms in formula \cref{eq:fracderivative1}, which may lead to tedious and time consuming calculation; nevertheless, from the definition of Gr\"{u}mwald-Letnikov fractional order derivative, the contribution of $f(x)$ near the initial point $x=0$ can be ignored if choosing an appropriate memory length in the approximation formula. In other words, for any $x\in (0,T]$, $f(x)$ is going to be taken into account only in the interval $[x-a,x]$, where $a\in (0,x)$ is called the memory length. Then, for an appropriate choice of $a$, we have
\begin{equation*}
  {}_0D^{\alpha}_x f(x) \approx {}_{x-a}D^{\alpha}_x f(x), \quad x>a,
\end{equation*}
and the following error estimate holds \cite{podlubny1997numerical}:
\begin{equation}\label{eq:shortmemoryerror1}
  \Delta (x): = |{}_0D^{\alpha}_x f(x) - {}_{x-a}D^{\alpha}_x f(x)| \leq\frac{Ma^{-\alpha}}{|\Gamma(1-\alpha)|},
\end{equation}
where $0<\alpha\neq0,1,2,...,$ and $M=\sup\limits_{0\leq x\leq L}|f(x)|$.

Similarly, for any $x\in [0,T)$  we have
\begin{equation*}
  {}_xD^{\alpha}_L f(x) \approx {}_{x}D^{\alpha}_{x+a} f(x), \quad x<L-a.
\end{equation*}

According to \cref{eq:shortmemoryerror1}, an appropriate memory length $a$ can be chosen to obtain a desired accuracy of the approximation formula.

By \cref{eq:fracderivative1}, \cref{eq:fracderivative1-R} and the Short Memory Principle, we introduce the following discrete approximations of fractional derivatives ${}_0D^{\alpha}_x$ and ${}_xD^{\alpha}_L$:
\begin{align*}
 {}_0\mathscr{D}^{\alpha}_x := & h^{-\alpha}\sum^{N-1}_{k=0}\omega_{k}^{(\alpha)}\Big\{f(x-kh) \\
                               & +\frac{1}{4}\alpha\Big(f(x-(k-1)h)-f(x-(k+1)h)\Big) \\
                               & +\frac{1}{8}\alpha^2\Big(f(x-(k-1)h)-2f(x-kh) \\
                               &\quad\quad\quad +f(x-(k+1)h)\Big)\Big\},\\
 {}_x\mathscr{D}^{\alpha}_L := & h^{-\alpha}\sum^{N-1}_{k=0}\omega_{k}^{(\alpha)}\Big\{f(x+kh) \\
                               & +\frac{1}{4}\alpha\Big(f(x+(k-1)h)-f(x+(k+1)h)\Big) \\
                               & +\frac{1}{8}\alpha^2\Big(f(x+(k-1)h)-2f(x+kh) \\
                               &\quad\quad\quad +f(x+(k+1)h)\Big)\Big\},
\end{align*}
where integer $N$ denotes the memory length, with
$$3\leq  N\leq \min\{n_1,n_2\}.$$
Then, for any $x,y\in (0, L)$, the discrete  two-sided fractional derivatives are expressed as
\begin{align*}
  \mathscr{D}^{\alpha}_x :=& \frac{1}{2}({}_0\mathscr{D}^{\alpha}_x+{}_x\mathscr{D}^{\alpha}_L)\\
  =& C_0h^{-\alpha}f(x) +h^{-\alpha} \sum^{N-2}_{j=1}C_j\Big(f(x-jh)+f(x+jh)\Big), \\
  \mathscr{D}^{\alpha}_y :=& \frac{1}{2}({}_0\mathscr{D}^{\alpha}_y+{}_y\mathscr{D}^{\alpha}_L)\\
  =& C_0h^{-\alpha}f(y) + h^{-\alpha}\sum^{N-2}_{j=1}C_j\Big(f(y-jh)+f(y+jh)\Big).
\end{align*}
where
{\footnotesize
\begin{align*}
C_0  = &   1 - \frac{\alpha^2}{2} - \frac{\alpha^3}{8},\\
C_1=&  \frac{\alpha}{8} + \frac{\alpha^2}{16}+\frac{1}{2\Gamma(-\alpha)}\Big(\frac{\Gamma(2-\alpha)}{2} (\frac{\alpha}{4}+\frac{\alpha^2}{8} ) \\
         &+{\Gamma(1-\alpha)} (1-\frac{\alpha^2}{4} ) +{\Gamma(-\alpha)} (-\frac{\alpha}{4}+\frac{\alpha^2}{8} )\Big),     \\
C_j  = &\frac{1}{2\Gamma(-\alpha)}\Big(\frac{\Gamma(j-\alpha+1)}{(j+1)!} (\frac{\alpha}{4}+\frac{\alpha^2}{8} ) \\
         &  +\frac{\Gamma(j-\alpha)}{j!} (1-\frac{\alpha^2}{4} ) +\frac{\Gamma(j-\alpha-1)}{(j-1)!} (-\frac{\alpha}{4}+\frac{\alpha^2}{8} )\Big)\\
&  \qquad \qquad  \qquad \qquad\qquad \qquad \qquad \qquad\text{for $j=2,3,\cdots, N-4$},\\
C_{N-3} =& \frac{\Gamma(N-\alpha-3)}{2\Gamma(-\alpha)(N-3)!}\Big(1-\frac{\alpha^2}{4}\Big)  + \frac{\Gamma(N-\alpha-4)}{2\Gamma(-\alpha)(N-4)!}\Big(-\frac{\alpha}{4}+\frac{\alpha^2}{8}\Big) , \\
C_{N-2} =& \frac{ \Gamma(N-\alpha-3)}{2\Gamma(-\alpha)(N-3)!}\Big(-\frac{\alpha}{4}+\frac{\alpha^2}{8}\Big).
\end{align*}
}

Correspondingly, from \cref{eq:definition-div-grad} we define the  discrete gradient $\widetilde{\nabla}^{\alpha}$ and the discrete divergence $\widetilde{\text{div}}^{\alpha}$  respectively as
\begin{align}\label{discrete-grd-div-1}
\widetilde{\nabla}^{\alpha}u :=& (\mathscr{D}^{\alpha}_x u, \mathscr{D}^{\alpha}_y u), \\
\widetilde{\text{div}}^{\alpha}\vec{f} :=& \widetilde{\nabla}^{\alpha}\cdot\vec{f} = \mathscr{D}^{\alpha}_x f_1 + \mathscr{D}^{\alpha}_y f_2.
\label{discrete-grd-div-2}
\end{align}

In what follows, we define the following  grids for the spatial domain $\Omega=[0,L]\times [0,L]$ and the temporal interval $[0,T]$. Given two positive integers $N_L$ and $N_T$, let $h=L/N_L$ and $\Delta t=T/N_T$ be the spatial and temporal step sizes, and set $x_k=y_k=k h$ for $k=0,1,\cdots, N_L$; $t_n=n\Delta t$ for $n=0,1,\cdots, N_T$.

Based on the  discrete approximations of the spatial fractional operators, \cref{discrete-grd-div-1}, \cref{discrete-grd-div-2} and the forward Euler scheme for the time discretization, we obtain the following fully discrete schemes of the proposed model.
\begin{discretization}
For $n = 0,1,2,...,N_T-1$ and $i,j=0,1,\cdots, N_L$,
\begin{align} \label{eq:dis1}
  \frac{u^{n+1}_{i,j}-u^n_{i,j}}{\Delta t} &= -~\widetilde{\text{div}}^{\alpha}(\eta^n_{i,j}),
\end{align}
where
\begin{align*}
  \eta^n_{i,j} := g(|\widetilde{\nabla}^{\beta}u^n_{i,j} |)\widetilde{\nabla}^{\alpha}u^n_{i,j} ,
\end{align*}
and $u^n_{i,j}$ denotes the approximation of $u(x_i,y_j,t_n)$ for any $i,j,n$.
\end{discretization}

\section{Experimental results}\label{section:Experimental results}
In this section, several numerical experiments are conducted to verify the effectiveness of the model. The method will be compared with some classic methods and fractional methods, such as the Gaussian filter, the Median filter, the classic PM model\cite{perona1990scale}, the classic TV(ROF) model \cite{rudin1992nonlinear}, the directional diffusion (DD) model and the fractional Fourier (FF) method \cite{yu2018image}.

There are many benchmark indexes to evaluate the quality of a processed image with the original one. Here  we will use the Peak Signal to Noise Ratio (PSNR) and the Structural Similarity index (SSIM) to assess the performance of the proposed method. The PSNR is an objective metric to estimate the fidelity between processed image and original image, and is easy to calculate with
\begin{equation*}
  PSNR(dB) = 10log_{10}\Big[\frac{(2^n-1)^2}{MSE}\Big],
\end{equation*}
where $n$ represents the bit of the image, usually chosen as $n=8$, and $MSE$ is the Mean Square Error between the processed image $u$ and the original one $u^\star$, given by
\begin{equation*}
  MSE = \frac{1}{N_L^2}\sum^{N_L}_{i=1}\sum^{N_L}_{j=1}\Big(u(i,j)-u^\star(i,j)\Big)^2.
\end{equation*}
It is easy to see that the fidelity of the processed image $u$   magnifies with PSNR increasing. It is natural that a processed image possesses an outstanding visual perception effect if it has almost no distortion compared with the original version, that is, it has a large PSNR.

Human visual system is highly sensitive about image structure, which indicates that the SSIM is another appropriate index to evaluate the performance of the method. From the perspective of image composition, the structure of image is defined as an attribution, which is independent of brightness and contrast of the image, to reflect the structure of objects. Such a definition    is extensively applied to evaluate the quality of videos and pictures.The SSIM, a metric to measure the similarity between two images, is described as
\begin{equation*}
  SSIM = \frac{(2\mu_u\mu_{u^\star}+c_1)(2\sigma_{uu^\star}+c_2)}{(\mu^2_u+\mu^2_{u^\star}+c_1)(\sigma^2_u+\sigma^2_{u^\star}+c_2)},
\end{equation*}
where $c_1 = (k_1l)^2$, $c_2 = (k_2l)^2$ with  $k_1 = 0.01$, $k_2 = 0.03$,   $l$ is the dynamic range of pixel values in image, $\mu_u$ and$\mu_{u^\star}$  are respectively the mean values of $u$ and $u^\star$, $\sigma_{uu^\star}$ is the covariance between $u$ and $u^\star$, and $\sigma^2_u$ and $\sigma^2_{u^\star}$  are respectively the variances of $u$ and $u^\star$. If the SSIM between $u$ and $u^\star$ is large, then $u$ is said to possess high similarity with respect to $u^\star$, which means that the method used to denoise is effective on structure preserving.

In terms of the estimate \cref{eq:shortmemoryerror1}, one needs to keep good balance between the memory length number $N$ and the  fractional order $\alpha$  so as to obtain an appropriate approximation of fractional derivative. In the model \cref{eq:origin}, the fractional order $\beta$ is used to control the edge-stopping function. The smaller the value of $\beta$, the less features the edge-stopping function will capture; and the larger the value, the more features the edge-stopping function will capture, but more noise will be captured at the same time, resulting in reduction of denoising performance. Hence, in all the numerical experiments with the scheme \cref{eq:dis1}, we take the spatial step size $h = 1$, the time step size $\Delta t = 0.5$, the memory length number $N = 15$, the fractional orders $\alpha = 1.67$ and $\beta = 1.55$, and the edge-stopping function $g(r)$ of the form \cref{eq:g1} with $\gamma=2$. With this setting, the proposed model is able to handle with most of cases well.

\begin{figure}[htbp!]
  \centering
  \subfloat[Lenna]{\includegraphics[width= 1.2in]{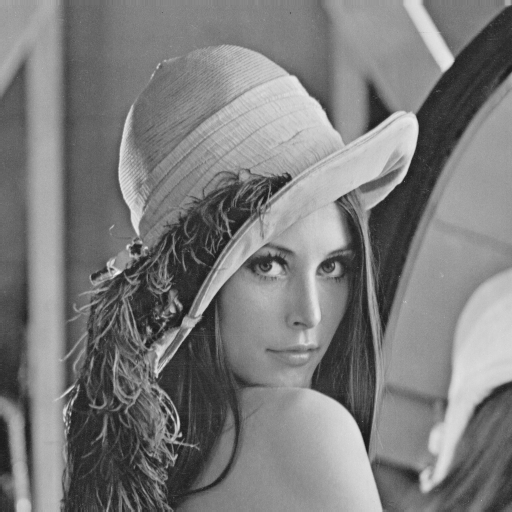}}\quad\quad
  \subfloat[Finger]{\includegraphics[width= 1.2in]{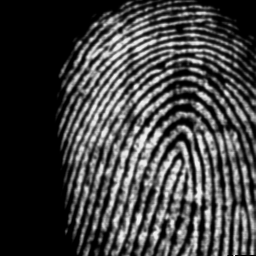}}\quad\quad
  \subfloat[Barbara]{\includegraphics[width= 1.4in]{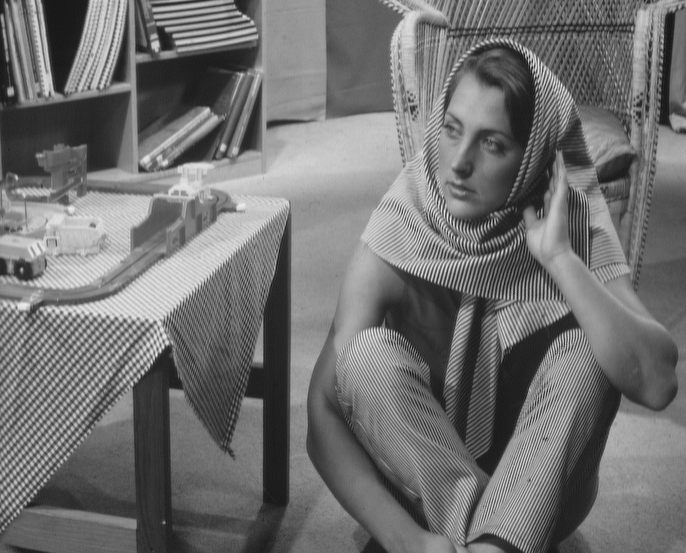}}\quad\quad
  \subfloat[Baboon]{\includegraphics[width= 1.2in]{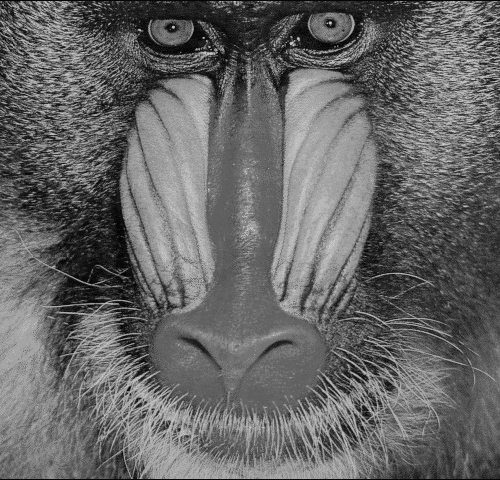}}\quad\quad
  \subfloat[Pepper]{\includegraphics[width= 1.2in]{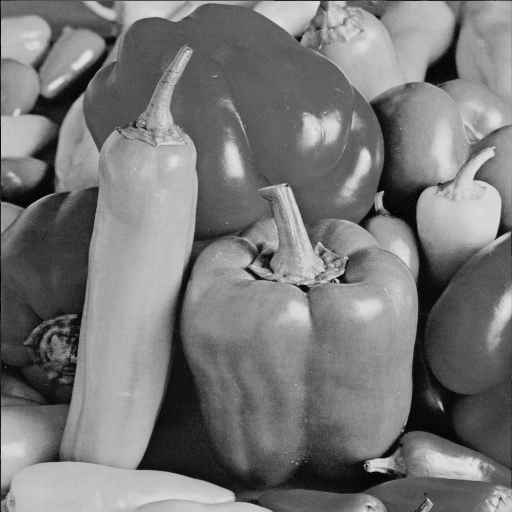}} \\
  \caption{Original standard test images.}\label{fig:original1}
\end{figure}

In \cref{fig:original1}, five standard test images are adopted to examine the proposed method. These pictures have large differences in feature, especially in edges and textures. The figures are corrupted with additive Gaussian white noise with standard deviations $\delta = 10,  15, 20, 25$.

\begin{figure}[htbp!]
  \centering
  \subfloat[Noise image]{\includegraphics[width= 1.5in]{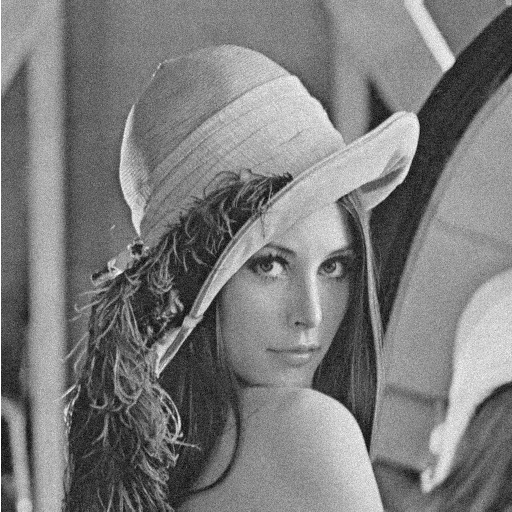}}\quad
  \subfloat[Gauss]{\includegraphics[width= 1.5in]{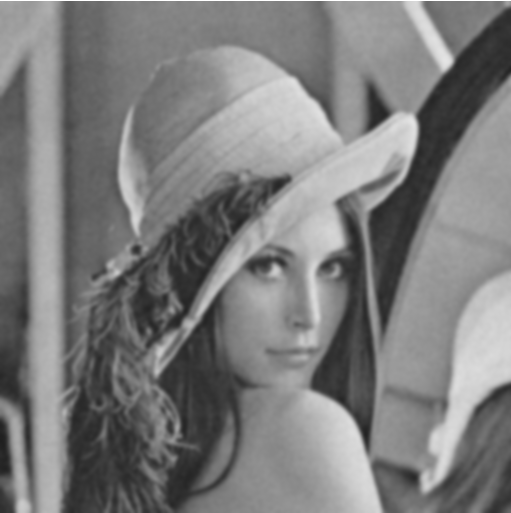}}\quad
  \subfloat[Median]{\includegraphics[width= 1.5in]{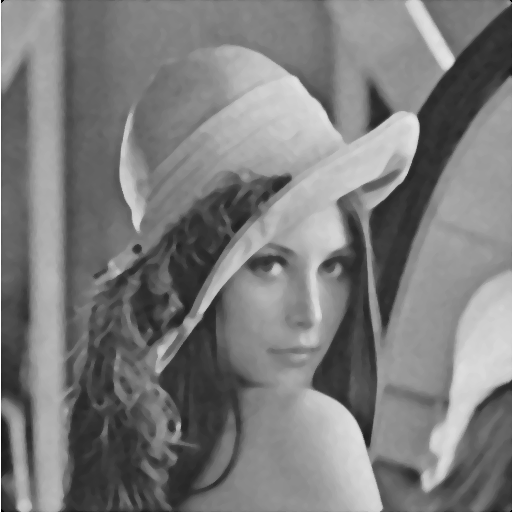}}\quad
  \subfloat[DD]{\includegraphics[width= 1.5in]{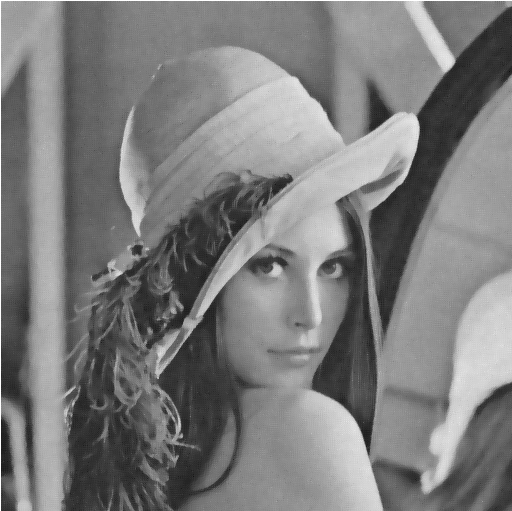}}\quad
  \subfloat[TV]{\includegraphics[width= 1.5in]{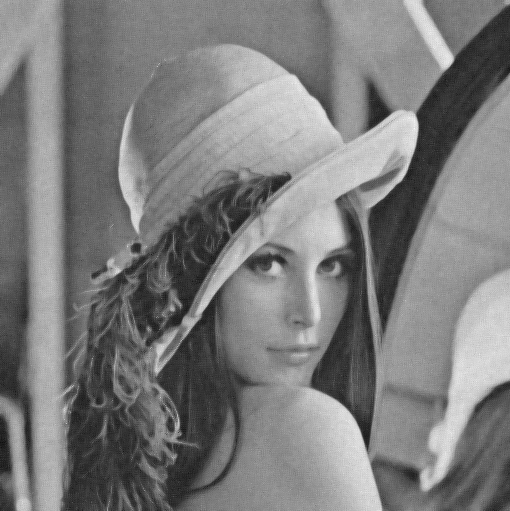}}\quad
  \subfloat[PM]{\includegraphics[width= 1.5in]{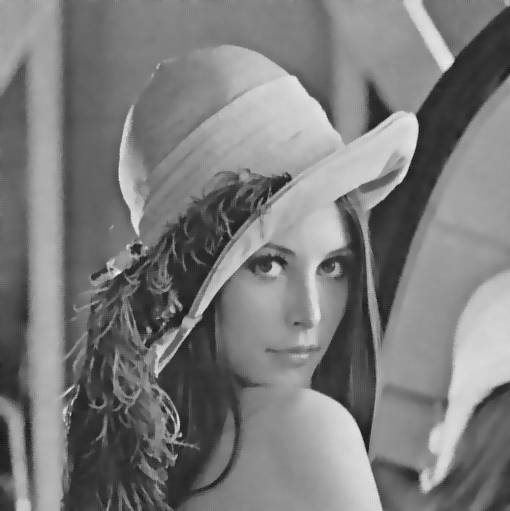}}\quad
  \subfloat[FF]{\includegraphics[width= 1.5in]{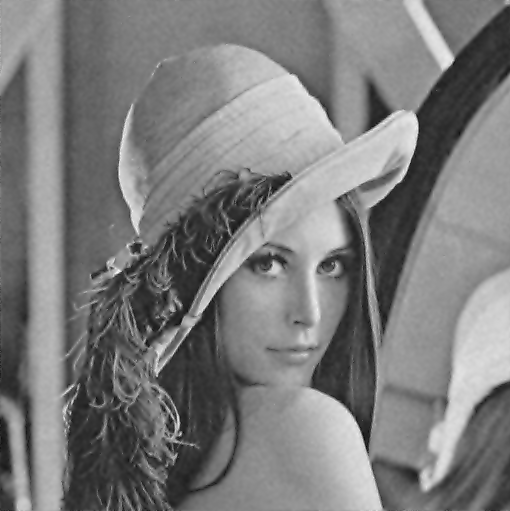}}\quad
  \subfloat[Proposed]{\includegraphics[width= 1.5in]{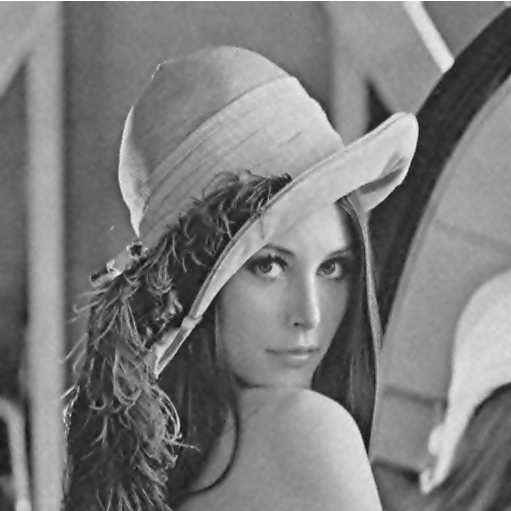}} \\
  \caption{Lenna image with additive Gaussian white noise $\delta=10$.}\label{fig:lennanoise1}
\end{figure}

\begin{figure}[htbp!]
  \centering
  \subfloat[Noise image]{\includegraphics[width= 1.5in]{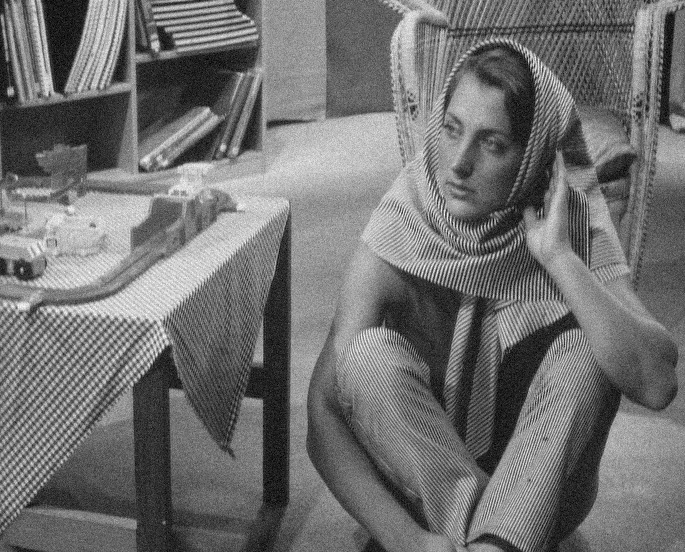}}\quad
  \subfloat[Gauss]{\includegraphics[width= 1.5in]{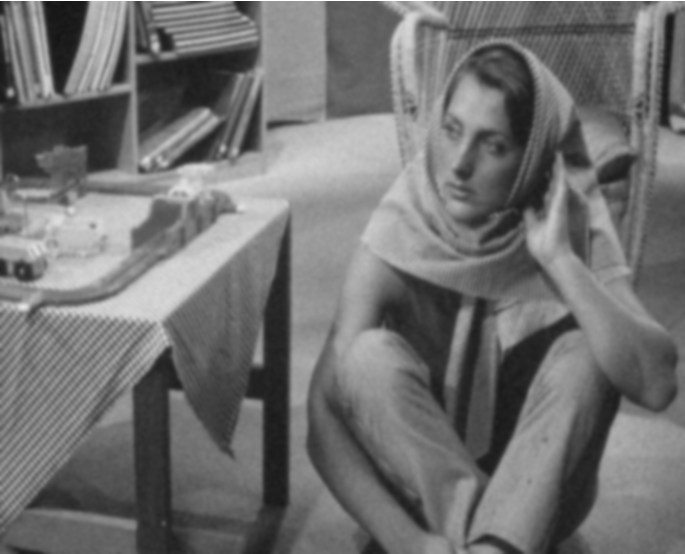}}\quad
  \subfloat[Median]{\includegraphics[width= 1.5in]{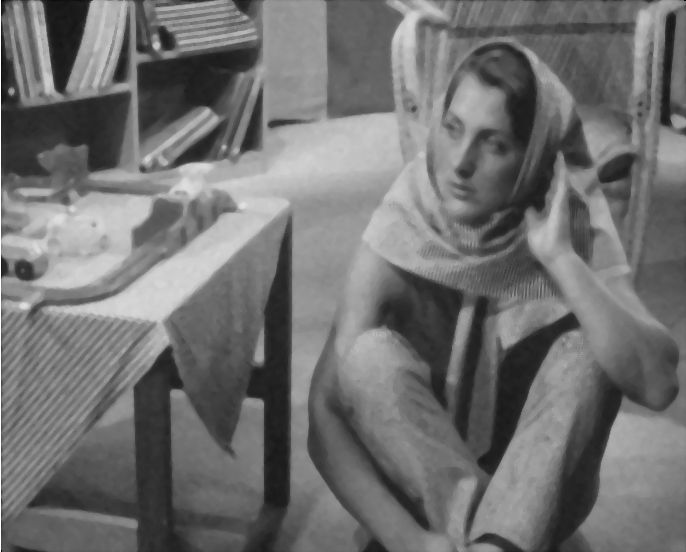}}\quad
  \subfloat[DD]{\includegraphics[width= 1.5in]{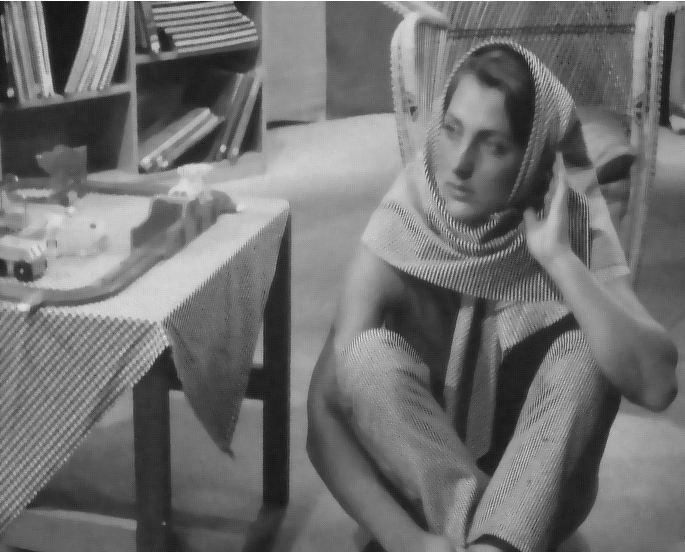}}\quad
  \subfloat[TV]{\includegraphics[width= 1.5in]{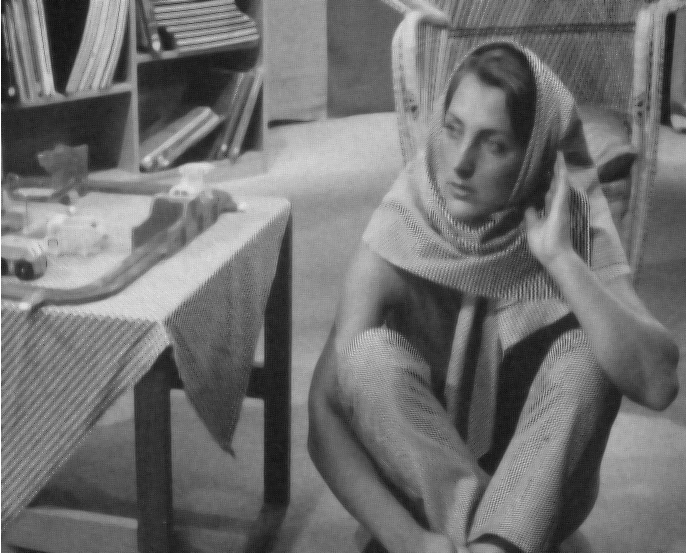}}\quad
  \subfloat[PM]{\includegraphics[width= 1.5in]{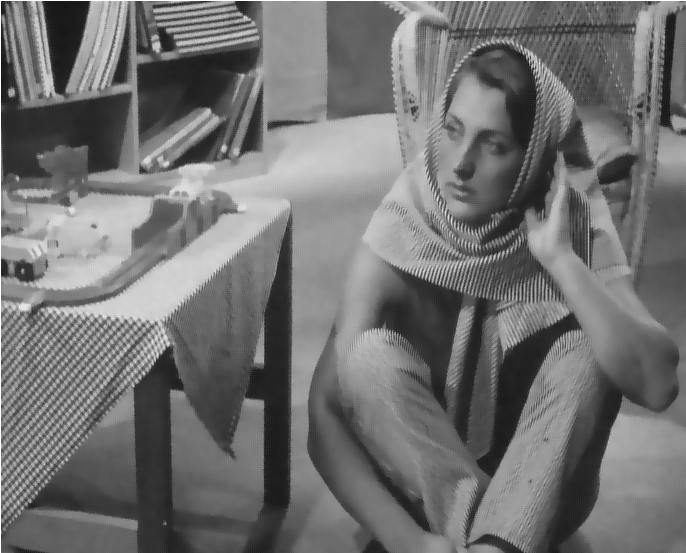}}\quad
  \subfloat[FF]{\includegraphics[width= 1.5in]{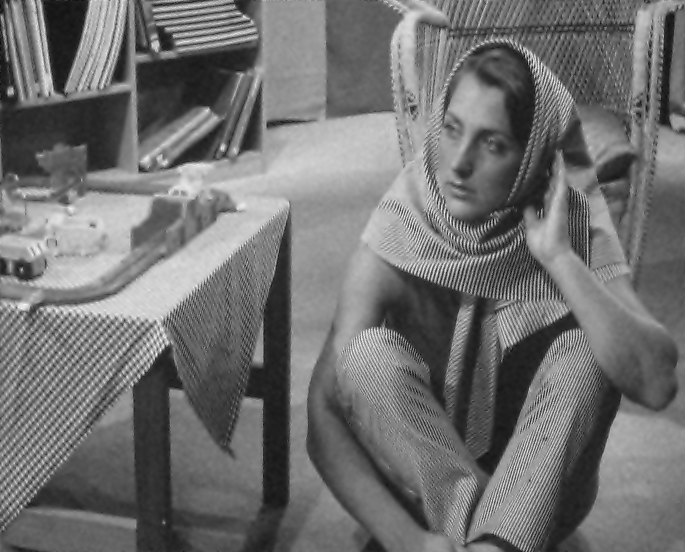}}\quad
  \subfloat[Proposed]{\includegraphics[width= 1.5in]{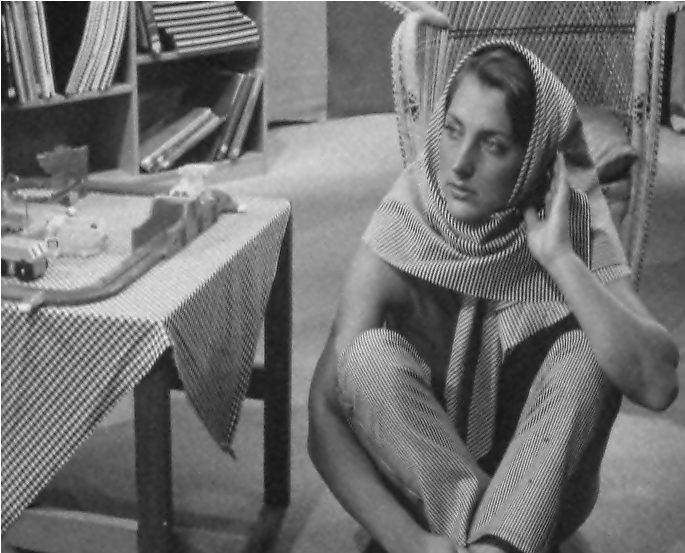}} \\
  \caption{Barbara image with additive Gaussian white noise $\delta=10$.}\label{fig:barbaranoise1}
\end{figure}

\begin{figure}[htbp!]
  \centering
  \subfloat[Noise image]{\includegraphics[width= 1.5in]{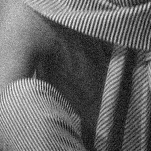}}\quad
  \subfloat[Gauss]{\includegraphics[width= 1.5in]{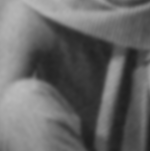}}\quad
  \subfloat[Median]{\includegraphics[width= 1.5in]{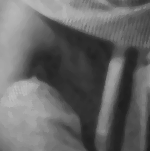}}\quad
  \subfloat[DD]{\includegraphics[width= 1.5in]{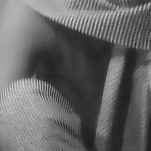}}\quad
  \subfloat[TV]{\includegraphics[width= 1.5in]{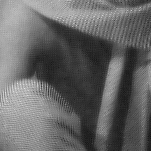}}\quad
  \subfloat[PM]{\includegraphics[width= 1.5in]{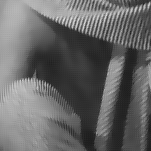}}\quad
  \subfloat[FF]{\includegraphics[width= 1.5in]{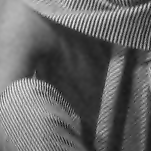}}\quad
  \subfloat[Proposed]{\includegraphics[width= 1.5in]{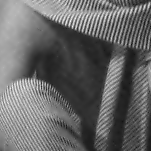}} \\
  \caption{Local barbara image with additive Gaussian white noise $\delta=10$.}\label{fig:barbaranoise2}
\end{figure}

\begin{figure}[htbp!]
  \centering
  \subfloat[Noise image]{\includegraphics[width= 1.5in]{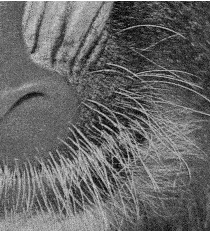}}\quad
  \subfloat[Gauss]{\includegraphics[width= 1.5in]{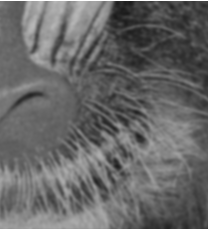}}\quad
  \subfloat[Median]{\includegraphics[width= 1.5in]{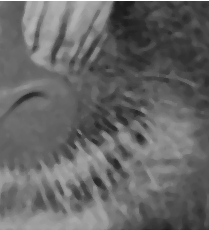}}\quad
  \subfloat[DD]{\includegraphics[width= 1.5in]{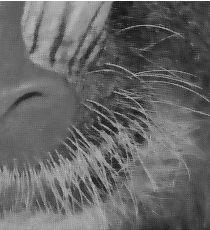}}\quad
  \subfloat[TV]{\includegraphics[width= 1.5in]{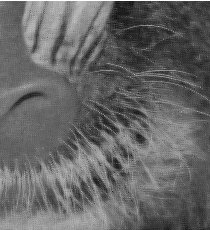}}\quad
  \subfloat[PM]{\includegraphics[width= 1.5in]{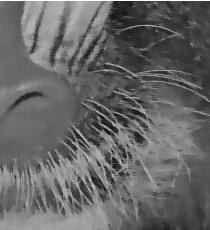}}\quad
  \subfloat[FF]{\includegraphics[width= 1.5in]{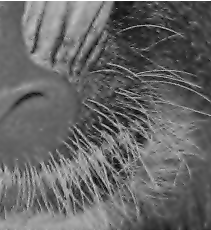}}\quad
  \subfloat[Proposed]{\includegraphics[width= 1.5in]{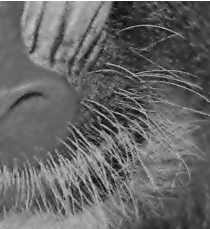}} \\
  \caption{Local baboon image with additive Gaussian white noise $\delta=10$.}\label{fig:baboonnoise1}
\end{figure}

With $\delta = 10$, \cref{fig:lennanoise1} and \cref{fig:barbaranoise1} show the denoising results of Lenna and Barbara images, respectively.\cref{fig:barbaranoise2} and \cref{fig:baboonnoise1} display two local parts of Barbara and Baboon images, respectively, with rich edges and textures.
From these figures we have the following observations.
 \begin{itemize}
 \item The Gaussian filter (\cref{fig:lennanoise1}(b) and \cref{fig:barbaranoise1}(b)) obviously blurs the edges and textures, for it actually is an isotropic filter which treats all the features and plain areas as same, e.g. the stripe textures on scarf and trousers (\cref{fig:barbaranoise2}(b)) are missing, and the moustaches of Baboon (\cref{fig:baboonnoise1}(b)) are adhesive. These two instances, \cref{fig:barbaranoise2}(b) and \cref{fig:baboonnoise1}(b), imply that the results of Gaussian filter with the isotropic diffusion process are not acceptable due to missing detailed information.

 \item The Median filter (\cref{fig:lennanoise1}(c) and \cref{fig:barbaranoise1}(c)) also blurs the edges and textures, but the burnishing effect of this filter is more temperate than that of the Gaussian filter. This situation happens for the nonlinear smoothing effect and anisotropy effect of the Median filter. The anisotropic diffusion process has an ability of feature preserving.  As shown in \cref{fig:barbaranoise2}(c), the stripe textures start to show up a little bit on the scarf, but are still missing  on trousers. Note that the result of the Median filter (\cref{fig:baboonnoise1}(c)) is even worse than  that of the Gaussian filter (\cref{fig:baboonnoise1}(b), since the anisotropic diffusion process of the Median filter has no feature selectivity for feature preserving at different pixels.

 \item Compared with the Gaussian filter and Median filter, the directional diffusion (DD) model has different results, as shown in \cref{fig:lennanoise1}(d) and \cref{fig:barbaranoise1}(d), with more edges and textures retained. The hair of Lenna has more hairlines in \cref{fig:lennanoise1}(d), while it is indistinct in \cref{fig:lennanoise1}(b) and \cref{fig:lennanoise1}(c). The same situation happens to the stripe textures on scarf and trousers of Barbara picture (\cref{fig:barbaranoise1}(d) and \cref{fig:barbaranoise2}(d)). Different from the Median filter, the anisotropic diffusion process of the DD model is feature-selective, because the diffusion process occurs only along the direction of tangent of isolux lines. The moustaches of Baboon in \cref{fig:baboonnoise1}(d) also implies that the DD model can preserve features due to its selective anisotropic diffusion effect.

 \item The TV model lightly alleviates the blurring effect and results are acceptable but not so satisfying. Benefited from $L_1$ regularization and total variation, the TV model is capable of retaining features. \cref{fig:lennanoise1}(e) and \cref{fig:barbaranoise1}(e) show that some details, e.g. the hair of Lenna and the stripe textures on scarf and trousers of Barbara, are preserved. However, there are still some noise remained on \cref{fig:barbaranoise2}(e) and \cref{fig:baboonnoise1}(e). In fact, as an isotropic model, the TV model is not easy to hold a superior balance of noise removal and feature preserving, and choosing an appropriate weight's setting between the total variation term and the fidelity term is usually empirical and difficult.

 \item The classic PM model has a good performance on edge preserving, but still obscures the textures and suffers from the staircase effect (\cref{fig:lennanoise1}(f)). In \cref{fig:barbaranoise1}(f), the stripe texture is retained well on the scarf through the classic PM method, whereas a portion of the texture on the trousers is still corrupted. \cref{fig:barbaranoise2}(f) shows that there are more stripe textures on the scarf and trousers, but the stripe textures on necktie are still missing. In \cref{fig:baboonnoise1}(f), there is almost no noise, but the staircase effect happens, e.g. on the nose of Baboon. These results imply that the classic PM model has an ability to retain features, yet it is still not excellent; nevertheless, compared with the Gaussian filter, the Median filter and the TV model, the classic PM model behaves better because of its anisotropy.

 \item The fractional Fourier (FF) method preserves edges and textures well in denoising process. The hairlines of Lenna in \cref{fig:lennanoise1}(g), the stripe textures on the scarf and trousers of Barbara in \cref{fig:barbaranoise1}(g), and the moustaches of Baboon in \cref{fig:baboonnoise1}(g) are clear. And there are enough stripe textures on the scarf and trousers remained in \cref{fig:barbaranoise2}(g), although some stripe textures on the necktie are lost.

 \item Our method keeps a more remarkable balance between noise removal and feature preserving, which suggests there is little damage on features when denoising; see \cref{fig:lennanoise1}(h) and \cref{fig:barbaranoise1}(h). Benefited from the using of two-sided fractional order derivatives, the proposed method also has a stronger ability of feature capturing than the  FF method. According to \cref{fig:lennanoise1}(h) and \cref{fig:barbaranoise1}(h), there are almost no staircase effect. In \cref{fig:barbaranoise2}(h) the stripe textures on the scarf and trousers of Barbara are remained, and our method has more stripe textures on the necktie than the FF method. \cref{fig:baboonnoise1}(h) also shows our method is superior, for the moustaches of Baboon are not adhesive. Metric indexes in \cref{table:experimentalresults1} to \cref{table:experimentalresults4} will also verify the superiority of our method.
 \end{itemize}

\begin{table}[htbp!]
  \caption{Experimental results with additive Gaussian white noise $\delta=10$.}
  \label{table:experimentalresults1}
  \centering
  \scalebox{0.75}{
  \begin{tabular}{|c |c |c |c |c | c | c | c | c | c | c| c| c| c| c| c|}
  \hline {Image} & {Index} &{Gauss} & {Median} & {DD} & {TV} & {PM} & {FF} & {Proposed} \\
  \hline
  \multirow{2}{*}{Lenna}
  & {PSNR} & {28.44} & {29.04} & {31.40} & {30.23} & {32.68} & {33.38} & {$\mathbf{34.01}$} \\
  & {SSIM} & {0.8120} & {0.8104} & {0.8125} & {0.7633} & {0.8477} & {0.8748} & {$\mathbf{0.8855}$} \\
  \hline
  \multirow{2}{*}{Finger}
  & {PSNR} & {21.61} & {20.26} & {26.63} & {25.15} & {30.05} & {30.07} & {$\mathbf{30.72}$} \\
  & {SSIM} & {0.6764} & {0.7819} & {0.7592} & {0.7413} & {0.7819} & {0.7861} & {$\mathbf{0.7958}$} \\
  \hline
  \multirow{2}{*}{Barbara}
  & {PSNR} & {24.18} & {24.08} & {26.68} & {25.13} & {27.35} & {30.21} & {$\mathbf{31.52}$} \\
  & {SSIM} & {0.6930} & {0.6821} & {0.7622} & {0.7071} & {0.7906} & {0.8740} & {$\mathbf{0.8826}$} \\
  \hline
  \multirow{2}{*}{Baboon}
  & {PSNR} & {22.01} & {21.06} & {25.97} & {23.45} & {26.01} & {26.41} & {$\mathbf{29.37}$} \\
  & {SSIM} & {0.5350} & {0.4342} & {0.7483} & {0.6799} & {0.7292} & {0.7535} & {$\mathbf{0.8775}$} \\
  \hline
  \multirow{2}{*}{Pepper}
  & {PSNR} & {29.73} & {30.51} & {32.20} & {30.78} & {33.22} & {33.44} & {$\mathbf{33.70}$} \\
  & {SSIM} & {0.8317} & {0.8241} & {0.8203} & {0.7778} & {0.8463} & {0.8550} & {$\mathbf{0.8593}$} \\
  \hline
  \end{tabular}
  }
\end{table}

\begin{table}[htbp!]
  \caption{Experimental results with additive Gaussian white noise $\delta=15$.}
  \label{table:experimentalresults2}
  \centering
  \scalebox{0.75}{
  \begin{tabular}{|c |c |c |c |c | c | c | c | c | c | c| c| c| c| c| c|}
  \hline {Image} & {Index} &{Gauss} & {Median} & {DD} & {TV} & {PM} & {FF} & {Proposed} \\
  \hline
  \multirow{2}{*}{Lenna}
  & {PSNR} & {28.31} & {28.64} & {29.73} & {28.49} & {31.24} & {31.29} & {$\mathbf{32.12}$} \\
  & {SSIM} & {0.8036} & {0.7944} & {0.7205} & {0.6733} & {0.8073} & {0.8205} & {$\mathbf{0.8522}$} \\
  \hline
  \multirow{2}{*}{Finger}
  & {PSNR} & {18.13} & {18.77} & {25.30} & {21.10} & {27.18} & {26.98} & {$\mathbf{27.84}$} \\
  & {SSIM} & {0.4698} & {0.7033} & {0.7059} & {0.6170} & {0.7292} & {0.7254} & {$\mathbf{0.7378}$} \\
  \hline
  \multirow{2}{*}{Barbara}
  & {PSNR} & {23.71} & {23.59} & {26.51} & {24.97} & {26.82} & {28.11} & {$\mathbf{28.52}$} \\
  & {SSIM} & {0.6594} & {0.6483} & {0.7262} & {0.6723} & {0.7620} & {0.8209} & {$\mathbf{0.8364}$} \\
  \hline
  \multirow{2}{*}{Baboon}
  & {PSNR} & {21.03} & {20.60} & {25.40} & {22.98} & {25.47} & {25.62} & {$\mathbf{27.86}$} \\
  & {SSIM} & {0.4224} & {0.3765} & {0.7249} & {0.6302} & {0.7177} & {0.7424} & {$\mathbf{0.7801}$} \\
  \hline
  \multirow{2}{*}{Pepper}
  & {PSNR} & {28.78} & {28.93} & {30.20} & {28.90} & {31.61} & {31.97} & {$\mathbf{32.02}$} \\
  & {SSIM} & {0.7905} & {0.7984} & {0.7339} & {0.6728} & {0.8083} & {0.8136} & {$\mathbf{0.8301}$} \\
  \hline
  \end{tabular}
  }
\end{table}

\begin{table}[htbp!]
  \caption{Experimental results with additive Gaussian white noise $\delta=20$.}
  \label{table:experimentalresults3}
  \centering
  \scalebox{0.75}{
  \begin{tabular}{|c |c |c |c |c | c | c | c | c | c | c| c| c| c| c| c|}
  \hline {Image} & {Index} &{Gauss} & {Median} & {DD} & {TV} & {PM} & {FF} & {Proposed} \\
  \hline
  \multirow{2}{*}{Lenna}
  & {PSNR} & {27.64} & {27.65} & {28.38} & {27.27} & {29.93} & {29.96} & {$\mathbf{30.75}$} \\
  & {SSIM} & {0.7835} & {0.7699} & {0.6670} & {0.6048} & {0.7701} & {0.8076} & {$\mathbf{0.8154}$} \\
  \hline
  \multirow{2}{*}{Finger}
  & {PSNR} & {17.98} & {17.87} & {23.07} & {19.64} & {25.08} & {25.11} & {$\mathbf{26.10}$} \\
  & {SSIM} & {0.4423} & {0.5922} & {0.6487} & {0.5425} & {0.6880} & {0.6887} & {$\mathbf{0.7031}$} \\
  \hline
  \multirow{2}{*}{Barbara}
  & {PSNR} & {23.37} & {23.18} & {25.31} & {23.75} & {25.61} & {26.19} & {$\mathbf{27.15}$} \\
  & {SSIM} & {0.6347} & {0.6215} & {0.6511} & {0.5611} & {0.7090} & {0.7810} & {$\mathbf{0.7938}$} \\
  \hline
  \multirow{2}{*}{Baboon}
  & {PSNR} & {20.84} & {20.44} & {23.70} & {21.74} & {23.68} & {23.75} & {$\mathbf{25.61}$} \\
  & {SSIM} & {0.3910} & {0.3568} & {0.6299} & {0.5151} & {0.6342} & {0.6687} & {$\mathbf{0.7087}$} \\
  \hline
  \multirow{2}{*}{Pepper}
  & {PSNR} & {28.09} & {28.11} & {28.89} & {28.59} & {30.39} & {30.51} & {$\mathbf{30.67}$} \\
  & {SSIM} & {0.7725} & {0.7818} & {0.6862} & {0.6261} & {0.7763} & {0.7832} & {$\mathbf{0.7991}$} \\
  \hline
  \end{tabular}
  }
\end{table}

\begin{table}[htbp!]
  \caption{Experimental results with additive Gaussian white noise $\delta=25$.}
  \label{table:experimentalresults4}
  \centering
  \scalebox{0.75}{
  \begin{tabular}{|c |c |c |c |c | c | c | c | c | c | c| c| c| c| c| c|}
  \hline {Image} & {Index} &{Gauss} & {Median} & {DD} & {TV} & {PM} & {FF} & {Proposed} \\
  \hline
  \multirow{2}{*}{Lenna}
  & {PSNR} & {27.07} & {26.91} & {27.46} & {26.48} & {28.97} & {28.99} & {$\mathbf{29.98}$} \\
  & {SSIM} & {0.7665} & {0.7504} & {0.6658} & {0.6090} & {0.7456} & {0.7861} & {$\mathbf{0.8105}$} \\
  \hline
  \multirow{2}{*}{Finger}
  & {PSNR} & {17.92} & {17.49} & {20.88} & {19.25} & {23.57} & {23.59} & {$\mathbf{24.71}$} \\
  & {SSIM} & {0.4336} & {0.5514} & {0.5768} & {0.5168} & {0.6517} & {0.6658} & {$\mathbf{0.6759}$} \\
  \hline
  \multirow{2}{*}{Barbara}
  & {PSNR} & {23.31} & {23.03} & {23.64} & {22.98} & {24.62} & {25.52} & {$\mathbf{26.12}$} \\
  & {SSIM} & {0.6285} & {0.6087} & {0.5808} & {0.5141} & {0.6608} & {0.7325} & {$\mathbf{0.7529}$} \\
  \hline
  \multirow{2}{*}{Baboon}
  & {PSNR} & {20.71} & {20.63} & {22.11} & {20.88} & {22.56} & {23.17} & {$\mathbf{23.84}$} \\
  & {SSIM} & {0.3883} & {0.3808} & {0.5066} & {0.4176} & {0.5231} & {0.6171} & {$\mathbf{0.6750}$} \\
  \hline
  \multirow{2}{*}{Pepper}
  & {PSNR} & {27.41} & {28.02} & {28.04} & {26.82} & {29.45} & {29.50} & {$\mathbf{29.81}$} \\
  & {SSIM} & {0.7556} & {0.7721} & {0.6683} & {0.6125} & {0.7521} & {0.7739} & {$\mathbf{0.7972}$} \\
  \hline
  \end{tabular}
  }
\end{table}

\cref{table:experimentalresults1} to \cref{table:experimentalresults4} give experimental results of the indexes PSNR and SSIM at different standard deviations $\delta = 10,  15, 20, 25$, with the highest scores being in bold type. We can see that all the the highest scores belong to the proposed method, especially in images with ample features, for instance, the Barbara and Baboon images. In these two pictures, the PSNR of the proposed method exceeds others to indicate that the processed image by the method has more fidelity. The SSIM of the method significantly surpasses other models on the test images, which implies that the proposed model has extremely high structural retention property, i.e. the processed picture has a preferable visual effect and is more similar to the original one in the structure sense.

Take \cref{table:experimentalresults4} as an example. The original images are corrupted with the additive Gaussian white noise $\delta=25$. From the Lenna image, the classic PM method has   better results of PSNR and SSIM than the first four models, i.e. the Gaussian filter, the Median filter, the directional diffusion model and  the classic TV model. The fractional Fourier method works slightly better than the PM model. Among these methods, the proposed method owns the best results of PSNR and SSIM, which verifies the capability of feature preserving.
Pictures such as the Barbara and Baboon images with rich, multi-layered edges and textures have a high probability of losing detailed information when processed with an inappropriate method. For the Baboon image, our method has 23.84dB PSNR and 0.6750 SSIM, which are superior to others.

In addition, we compare our method with the works in \cite{jia2017new, litvinov2011modified} based on the TV-Stokes model in terms of PSNR and SSIM. Images are blurred by Additive Gaussian White Noise with $\delta=25$. \cref{table:experimentalresults5} shows that our method yields better results both in PSNR and SSIM.
\begin{table}[htbp!]
  \caption{Experimental results with additive Gaussian white noise $\delta=25$.}
  \label{table:experimentalresults5}
  \centering
  \scalebox{0.75}{
  \begin{tabular}{|c |c |c |c |c |}
  \hline {Image} & {Index} &{\cite{litvinov2011modified}} & {\cite{jia2017new}} & {Proposed} \\
  \hline
  \multirow{2}{*}{Lenna}
  & {PSNR} & {28.87} & {29.21} & {$\mathbf{29.98}$} \\
  & {SSIM} & {0.7752} & {0.7906} & {$\mathbf{0.8105}$} \\
  \hline
  \multirow{2}{*}{Finger}
  & {PSNR} & {24.09} & {24.23} & {$\mathbf{24.71}$} \\
  & {SSIM} & {0.6551} & {0.6609} & {$\mathbf{0.6759}$} \\
  \hline
  \multirow{2}{*}{Barbara}
  & {PSNR} & {24.98} & {25.39} & {$\mathbf{26.12}$} \\
  & {SSIM} & {0.7197} & {0.7315} & {$\mathbf{0.7529}$} \\
  \hline
  \multirow{2}{*}{Baboon}
  & {PSNR} & {22.14} & {22.78} & {$\mathbf{23.84}$} \\
  & {SSIM} & {0.6547} & {0.6616} & {$\mathbf{0.6750}$} \\
  \hline
  \multirow{2}{*}{Pepper}
  & {PSNR} & {28.74} & {29.38} & {$\mathbf{29.81}$} \\
  & {SSIM} & {0.7595} & {0.7718} & {$\mathbf{0.7972}$} \\
  \hline
  \end{tabular}
  }
\end{table}

\section{Conclusion}\label{section:Conclusion}
The proposed PDE algorithm based on the two-sided spatial-fractional anisotropic diffusion equation has a satisfactory performance on feature preserving when denoising. According to the experimental results of some standard test images, the method possesses the best PSNR and SSIM at the same time, compared with other methods. The  main reason for the performance lies in the  following two aspects:
\begin{itemize}
\item Compared with integer order  derivatives, the fractional derivatives are more capable to nonlinearly preserve the low frequency information and enhance the high frequency information of images, thus keep better balance between noise removal and feature preserving;
\item Compared with one-sided fractional derivatives, the two-sided fractional derivatives capture more neighbourhood information, thus is more suitable to depict the local self-similarity of images.
\end{itemize}

\bibliographystyle{plain}
\bibliography{References}

\end{document}